\documentclass[11pt]{amsart}
\usepackage{amssymb}
\usepackage{amsfonts}
\usepackage{amsmath}
\usepackage[matrix,arrow,graph]{xy}

\theoremstyle{plain}
\newtheorem{thm}{Theorem}
\newtheorem{lemma}[thm]{Lemma}
\newtheorem{prop}[thm]{Proposition}
\newtheorem{cor}[thm]{Corollary}
\theoremstyle{definition}

\newtheorem*{defn}{Definition}

\theoremstyle{remark}
\newtheorem*{rem}{Remark}

\newcommand{\wt}{\widetilde}

\def\Span{\mathop{\rm span}}

\newcommand{\C}{{\mathbb C}}

\newcommand{\R}{{\mathbb R}}

\newcommand{\Q}{{\mathbb Q}}

\newcommand{\cal}{\mathcal}
\newcommand{\cS}{{\cal S}}

\newcommand{\cL}{{\cal L}}

\newcommand{\cM}{{\cal M}}
\newcommand{\IC}{{\mathbf {IC^{\textstyle \cdot}}}}
\newcommand{\mb}{\mathbf}

\begin{document}

\title{On the reducibility of characteristic varieties}
\author{Tom Braden}
\email{braden@math.harvard.edu}
\begin{abstract} We show that some monodromies in the Morse local systems of
a conically stratified perverse sheaf imply that other Morse local
systems for smaller strata do not vanish.  This result is then used
to explain the examples of reducible characteristic
varieties of Schubert varieties given by Kashiwara and Saito
in type $A$ and by Boe and Fu for the Lagrangian Grassmannian.
\end{abstract}
\maketitle

Let $X = G/P$ be a flag variety associated to a complex semisimple group $G$
and parabolic subgroup $P$.
Let $Y\subset X$ be a Schubert variety, the closure of an orbit under the
action of a Borel subgroup $B\subset G$ on $X$.  The characteristic cycle
$CC(\IC(Y))$ of the intersection cohomology sheaf on $Y$ is an object of
considerable representation theoretic interest.  It is an $n$-dimensional
closed-support cycle (where $n = \dim_\R X$) supported on the conormal variety 
$\Lambda \subset T^*X$ to the Schubert stratification.  Since 
$\Lambda$ is purely $n$-dimensional, with one component for
each stratum, this amounts to giving an integer for each Schubert cell.

The support of the characteristic cycle is known as the
characteristic variety.
In \cite{KL}, Kazhdan and Lusztig stated a conjecture 
equivalent to the statement that for $G=SL_n$, the
characteristic varieties of intersection cohomology sheaves should be 
irreducible.
This work arose from an attempt to understand the counterexample
found by Kashiwara and Saito \cite{KSa} to this conjecture.

We present a theorem (Theorem \ref{mainthm}, also see 
Corollary \ref{maincor}) which explains many of the examples
currently known (for general semisimple $G$, not just $G = SL_n$)
in which $CC(\IC(Y))$ is reducible.  It
gives a geometric condition on a stratified variety $Y$ and a stratum $S$
that implies that the component $[\Lambda_S]$ lying over $S$ appears in 
$CC(\IC(Y))$; it requires that $Y$ have a conical singularity
along $S$.

We use this result to verify Kashiwara and Saito's example, and
to explain the examples of non-reducible characteristic varieties
arising in Lagrangian Grassmannians found by Boe and Fu \cite{BF}.
Our result can also be used to explain all but one of the examples
found by Tanisaki \cite{T} of reducible characteristic
varieties in full flag manifolds.  
Our result cannot show that components of $\Lambda$ do
not occur in the characteristic variety, 
nor does it provide a way to calculate the 
actual multiplicities, as both Boe and Fu and Tanisaki were able to do.

We take $\Q$ as our ground field throughout this paper; all intersection
cohomology sheaves are taken with constant $\Q$ coefficients.
\section{Vanishing cycles of conical sheaves.}

This section presents an argument from \cite{thesis}.
Let $V$ be a complex vector space, and suppose an object ${\mb A}$ 
in the derived category $D^b(V)$
is constructible with respect to an algebraic, $\C^*$-conical
stratification.    
Let $f\colon V\to \C$ be a linear function, and 
let $W = f^{-1}(0)$.  For definitions of the sheaf-theoretic concepts we
use, including the vanishing cycles functor $\phi_f$ and the 
Fourier transform $F$, see \cite{KS}.

\begin{thm} \label{linvs} 
The multiplicity of $[T^*_{\{0\}}V]$ in $CC({\mb A})$
is the same as the multiplicity of $[T^*_{\{0\}}W]$ in $CC(\phi_f{\mb A})$.
\end{thm}
\begin{proof}  The Fourier transform $F\colon D^b(V) \to D^b(V^*)$ 
respects the characteristic cycle: we have 
$CC(F{\mb A}) = CC({\mb A})$, using the
natural identification $T^*V \cong V \times V^* \cong T^*V^*$ (this is given 
as Exercise 9.7 in \cite{KS}).  Since the
multiplicity of the zero section $T^*_{V^*}V^*$ in $CC(F{\mb A})$ is the 
Euler characteristic of the stalk cohomology of ${\mb A}$ at a point in 
the open stratum,
it will be enough to show that for generic inclusions $i_1$, $i_2$ 
of a point $\{p\}$ into $V^*$ and $W^*$, respectively, the restrictions 
$i_1^*F{\mb A}$ and $i_2^*F\phi_f{\mb A}$ are isomorphic.

We can consider $f$ as an element of $W^\bot \subset V^*$. Let
 $s\colon W \to W \times W^\bot$ be given by $s(w) = (w, f)$.  
Then $\phi_f$ is naturally isomorphic to following the top
row in the following diagram:

\[\xymatrix{D^b(V) \ar[r]^-{\nu_W}\ar[d]^F & D^b(W\times V/W)\ar[dr]^F
 \ar[r]^-F & 
D^b(W\times W^\bot) \ar[r]^-{s^*} \ar[d]^F & D^b(W)\ar[d]^F\\
D^b(V^*) \ar[rr]^-{\nu_{W^\bot}} & & D^b(W^*\times W^\bot)\ar[r]^-{(s')^*} & 
D^b(W^*)
}\]
Here $\nu_W$, $\nu_{W^\bot}$ are the specialization functors (using the 
natural identifications $T_WV \cong W \times V/W$ and $T_{W^\bot}\cong
W^\bot \times W^*$), the 
functors marked $F$ are the appropriate Fourier transforms, and
$s'(\omega) = (\omega, f)$.  The left quadrilateral and 
the triangle naturally commute by \cite{BG}, Propositions 2.3 and 2.4,
and the right hand square commutes by \cite{KS}, Proposition 3.7.13.

Thus $F\phi_f{\mb A} \cong (s')^*\nu_{W^\bot}F{\mb A}$.  It is easy now
to see that $\nu_{W^\bot}$ and $(s')^*$ both leave the dimensions 
of the stalk cohomology
at a generic point unchanged (for $(s')^*$, use the fact that
$\nu_{W^\bot}F{\mb A}$ is conical).
\end{proof}

Similar techniques have been useful in the calculation of certain 
categories of perverse sheaves \cite{B}, \cite{BG}.

\section{The main result}\label{mresult}

Suppose a smooth complex variety 
$V$ carries an algebraic Whitney stratification $\cS$.
Denote the conormal variety by $\Lambda \subset T^*V$, and let
$\Lambda_S = \overline{ T^*_SV}$ be the component of $\Lambda$ lying over 
$S \in \cS$.  Let $\wt\Lambda$ be the smooth part of $\Lambda$, and 
for each stratum $S$ put $\wt\Lambda_S = \wt\Lambda \cap \Lambda_S$.

\begin{defn}
  If two 
strata $S$ and $T$ satisfy $\dim_\C(\Lambda_S \cap \Lambda_T) =\dim_\C 
\Lambda -1$,
we say that $S$ and $T$ {\em meet microlocally in codimension one.}
\end{defn}

Given an $\cS$-constructible perverse sheaf ${\mb P}$ on $V$, we let
$\cM_S({\mb P})$ be the Morse or vanishing cycles local system
of ${\mb P}$ at the stratum $S$ (for a definition, see \cite{MV}).  
It is a local system on $\wt\Lambda_S$, whose
dimension is, up to a sign, the multiplicity of $[\Lambda_S]$ in 
$CC({\mb P})$.

Now suppose that $V$ is a complex vector space and the
stratification $\cS$ is $\C^*$-conical.  Let $S \ne \{0\}$ be a 
stratum, and take a point $x\in S$.  Define a loop $\gamma$ in
$\wt\Lambda_S$ by choosing a point $(x, \xi) \in \wt\Lambda_S$
and letting $\gamma(\theta) = (e^{ i\theta}, \xi)$, $\theta \in [0,2\pi]$.  
Note that any two such loops are homotopic.

\begin{thm} \label{mainthm} Suppose that ${\mb P}$ is an 
$\cS$-constructible perverse sheaf on $V$.  If 
$S$ meets $\{0\}$ microlocally in codimension one, and
the monodromy of $\cM_{S}({\mb P})$ 
around $\gamma$ is not multiplication by 
$(-1)^{d-1}$, then $\cM_0({\mb P}) \ne 0$.
\end{thm}

In the case ${\mb P} = \IC(\overline{S})$, we get:
\begin{cor} \label{maincor} 
Suppose $X \subset V$ is an even-dimensional conical
variety, $S$ is the smooth locus of $X$, and 
$\{0\}$ and $S$ meet microlocally in codimension one.  Then 
$\cM_0(\IC(X)) \ne 0$.
\end{cor}
\begin{proof} 
Just note that $\cM_S(\IC(\overline{S})) \cong \Q_{\wt\Lambda_S}$.  
\end{proof}

\begin{proof}[Proof of Theorem \ref{mainthm}]
Since we are free to choose the point $(x,\xi) \in \wt\Lambda_S$, we choose
it as follows. The dual cone $S^\vee \subset V^*\cong T^*_0V$ to $S$ is
defined to be $\Lambda_S \cap T^*_0V$; $S$ and $S^\vee$ are cones over
dual projective varieties in ${\mathbb P}V$ and ${\mathbb P}V^*$.
By assumption $S^\vee$ is a divisor and 
so it is not contained in any $T^\vee$ for any
stratum $T\ne S, \{0\}$.  Let  $\xi$ be
a smooth point of $S^\vee \setminus S^\vee \cap\bigcup_{R\ne S,\{0\}} R^\vee$.
Then  $L=\{v\in V\mid (v,\xi)\in \Lambda\}$
is a line through the origin contained in $S\cup\{0\}$; let $x$ be
any nonzero point in $L$. 

In microlocal language, the point $(0, \xi)$ represents a codimension
one point of $\Lambda$:  it is a smooth point of both $\Lambda_0$ and
$\Lambda_{S}$, they intersect transversely
there, and it lies in no other component of $\Lambda$. 

Let $f\colon V\to \C$ be the linear function with $df_x=\xi$.
Our choice of $\xi$ implies that $f|_S$ has singularities only
on $L$ and that if we choose a normal slice $N$ to $L$ at $x\in L$,
the singularity of $f|_{N\cap S}$ is Morse at $x$.

Let $L^\circ = L \setminus \{0\}$
and embed $L^\circ$ into $\wt\Lambda$
by the map $\alpha(y) = (y, df_y) = (y, \xi)$.  
Put ${\mb P}_f = \phi_f({\mb P})$; it is a conical perverse sheaf 
supported on $L$.  The local system on $L^\circ$
can be computed using stratified Morse theory:
\[{\mb P}_f|_{L^\circ} = \alpha^*\cM_{S}({\mb P}) \otimes \cL,\]
where $\cL = \phi_{(f|_S)}(\C_{S})$.  
But the local system $\cL$ is easy to compute by Picard-Lefchetz theory;
its monodromy is multiplication by 
$(-1)^{d-1}$, where $d = \dim_\C(S)$.  Thus the monodromy of 
${\mb P}_f|_{L^\circ}$ is nontrivial.

An elementary application of the theory of perverse sheaves on a complex
line (see \cite{V} and the example following the proof of Theorem 3.3 in 
\cite{MV}) implies that a perverse sheaf with nontrivial monodromy around the
origin must have a nonzero vanishing cycle at $0$; thus
$\cM_0({\mb P}_f) \ne 0$.
Then by Theorem \ref{linvs}, $\cM_0({\mb P}) = \cM_0({\mb P}_f)$, 
so we are done.
\end{proof}

\begin{rem} A strengthening of Theorem \ref{linvs} using 
Morse local systems instead of characteristic cycles gives
the slightly stronger conclusion that $\cM_0({\mb P})$ is a 
nontrivial local system.
\end{rem}

We have stated Theorem \ref{mainthm} assuming
that the ambient variety is a vector space and 
the smaller stratum is a point, but it can also be applied to a 
pair of strata $(S, T)$ in a general stratified complex manifold $V$ 
with $T \subset \overline{S}$, so long as the stratification is conical 
along $T$; simply intersect with a normal slice $N$ to $S$ at a point
$s\in S$, or, equivalently, apply the specialization functor 
$\nu_S\colon D^b(V) \to D^b(T_SV)$ and restrict to the fiber over $s$.

If $V\cong \C^n$ but the stratification of $V$ is not conical, one can still look at
the specialization ${\mb P}' =\nu_0({\mb P})$, which is conical. 
Since $\cM_0({\mb P}') \cong \cM_0({\mb P})$ (see \cite{Sa}), Theorem 2 may 
apply to ${\mb P}'$ to give information about ${\mb P}$.
We cannot use this to 
deduce that Theorem 2 holds without the conical assumption,
however.  There are stratifications for which 
$S$ meets $0$ microlocally in codimension 
one, but the smooth part of the specialization $\nu_0(\overline S)$
doesn't.  

For instance, if $\overline{S}$ is the 
variety $\{(x,y,z)\in \C^3\mid xy=z^3\}$, the specialization
${\mb P}' = \nu_S(\IC(\overline{S}))$ 
is supported on $\{xy = 0\}$.  A little calculation shows that 
the  composition series of the perverse sheaf ${\mb P}'$
consists of simple intersection cohomology sheaves with constant coefficients 
supported on 
$\{x = 0\}$, $\{y=0\}$, $\{x = y = 0\}$, and $\{(0,0,0)\}$.  
If Theorem \ref{mainthm} applied to ${\mb P}'$, it would
have to apply for one of these simple components, since 
the Morse local systems $\cM_S$ functors are exact on the 
category of perverse sheaves.  But this is clearly not
the case, since all but the last have $\cM_0 = 0$.  

This 
singularity appears in the flag variety for $G=G_2$, 
associated to the pair of Weyl group elements $(tst, tstst)$,
where $s$,$t$ are the reflections for the long and short simple roots,
respectively.

\section{Kashiwara and Saito's example}

In \cite{KSa} Kashiwara and Saito gave an example of a pair of Schubert
varieties $Z \subset Y$ in the variety of complete flags in $\C^8$
so that $\Lambda_Z$ appears in the support of $CC(\IC(Y))$.
We recall their description of a normal slice to $Y$ at a 
generic point of $Z$. Let $V\cong \C^{16}$ be the space  
${\rm Mat}(2\times 2, \C)^4$
of $4$-tuples of $2\times 2$ complex matrices.
Let $X$ be the variety
\[\{(A_0, A_1, A_2, A_3) \in V \mid \det(A_i) = 0 
\;\text{and}\; A_iA_{i+1} = 0\;\text{for all}\; i\,\},\]
where the index $i$ is taken modulo $4$.  
It is an $8$-dimensional conical subvariety of $V$.

We show that $X$ satisfies the hypotheses of Corollary \ref{maincor}.
Let $G = (GL_2(\C))^4 \times \C^*$ act on $V$ by
\[(g_0,g_1, g_2, g_3, t)\cdot(A_0,A_1,A_2,A_3) = 
(tg_0A_0g_1^{-1}, g_1A_1g_2^{-1},
g_2A_2g_3^{-1}, g_3A_3g_0^{-1}).\]
Then $X$ is $G$-invariant, and in fact is the closure of a 
$G$-orbit: $X = \overline{G\cdot(A,A,A,A)}$, where $A$ is any nonzero
nilpotent matrix.  The open orbit is the subset of points
$(A_0, A_1, A_2, A_3)$ in $X$ for which all the $A_i$ are nonzero.

\begin{prop} 
The smooth part of $X$ meets $\{0\}$ microlocally
in codimension one.  
\end{prop}
\begin{proof} 
This is equivalent to showing that the
 dual cone $X^\vee \subset V^*$ is a divisor.
Using the inner product on $V$ given by matrix coordinates, 
$V^*$ is naturally identified again with ${\rm Mat}(2\times 2, \C)^4$, with 
the action of $G$ given by 
\[(g_0, g_1, g_2, g_3, t)\cdot(A_0,A_1,A_2, A_3) =
(t^{-1}g_1A_0g_0^{-1}, g_2A_1g_1^{-1}, g_3A_2g_2^{-1}, g_0A_3g_3^{-1}).\]
Certainly $X^\vee$ is $G$-stable, and it is easily checked 
that $x = (I+A^t,I,I,I)$ lies in $X^\vee$, where again $A$ is any
nonzero nilpotent matrix.
Its stabilizer is $G_x = \{(g, g, g, g, 1)\mid g = aI + bA^t,\, a\ne 0\}$,
a two-dimensional group, so $\dim_\C X^\vee = 17 - 2 = 15$.
\end{proof}

Thus Corollary \ref{maincor} applies to this example, showing that
$[\Lambda_{\{0\}}]$ appears in $CC(\IC(X))$.  

Also note that $X$ is a toric variety: although the
maximal torus $T\subset SL_8(\C)$ which acts on the flag variety 
is only $7$-dimensional, there is a larger torus in $G$ which
acts on $X$ with finitely many orbits.  A laborious calculation with
an algorithm from \cite{thesis} or the equivalent formula
in \cite{GKZ} (Theorem 2.12 in Chapter 10) shows that 
$\dim\cM_0(\IC(X)) = 1$.  

\section{The Lagrangian Grassmannian}

In \cite{BF}, Boe and Fu computed the characteristic cycles $CC(\IC(Y))$
for Schubert varieties $Y$ in Hermitian symmetric spaces.  In the
compact cases, reducible characteristic varieties appeared only for
the Lagrangian Grassmannian $X$ of Lagrangian subspaces of 
a complex symplectic space $\C^{2n}$.  For concreteness, suppose
that the symplectic form is given by 
$\omega = \sum_{i=1}^n {\mb e}^*_i \wedge {\mb e}^*_{2n + 1 - i}$, where
the ${\mb e}^*_i$ form the dual basis to the standard basis of $\C^{2n}$.

The Schubert decomposition of $X$ is given as follows.
Let $B \subset {\rm Sp}(\C^{2n}, \omega)$ be the Borel group 
of transformations preserving the standard flag.
Given a word $w \in \{\alpha,\beta\}^n$ of length $n$ in the letters
$\alpha$ and $\beta$, we define a cell $S_w$ as follows.  Let
$\bar{w}\in \{\alpha,\beta\}^{2n} $ be the word for which
 $\bar{w}(2n + 1 -i)\ne \bar{w}(i) = w(i)$ for $1\le i \le n$.  
Then $E_w = \Span\{{\mb e}_i \mid
\bar{w}(i) = \alpha\}$ is a point in $X$, and we put $S_w = B\cdot E_w$ and
$X_w = \overline{S_w}$. 
Also let $N_w = B^-E_w$, where $B^-$ is the opposite Borel to $B$;  it is a  
normal slice to $S_w$ through $E_w$.

Let ${\mb P_w} = \IC(X_w)$.
Given $v, w \in  \{\alpha, \beta\}^n$, let $m_v^w$ be the multiplicity with which 
$[\Lambda_v]$ appears in $CC({\mb P_w})$.  Boe and Fu's calculation of these numbers
can be summarized as follows.
\begin{thm}[Boe and Fu \cite{BF}, Theorem 7.1D] \label{BFthm} 
$m_v^w = 0$ or $1$ for any $v,w \in \{\alpha, \beta\}^n$.  It is
$1$ if and only if there is a chain \[X_w = X_1 \supset X_2 \supset
\dots \supset X_n = X_v\] of Schubert varieties $X_i = X_{w_i}$ 
so that for all $1\le i < n$ the codimension 
$\dim_\C X_i - \dim_\C X_{i+1}$ is even and $S_{w_i}$ meets
$S_{w_{i+1}}$ microlocally in codimension one.
\end{thm}

Suppose that there exists a chain $\{w_i\}$ satisfying the conditions
of this theorem.  We will use Theorem \ref{mainthm} to give a shorter proof that 
$m_v^w \ne 0$.  Let $E_i = E_{w_i}$, $S_i = S_{w_i}$, $\wt\Lambda_i = \wt\Lambda_{S_i}$.
Assume inductively that $m^w_{w_i} > 0$ (this is trivially true for
$i=1$), so $\cM_i = \cM_{S_{i}}({\mb P_w})$ is a nonzero 
local system on $\wt{\Lambda}_i$.  We will show that this implies $\cM_{i+1}$
is nonzero as well.

Put $E = E_{i+1}$, $N = N_{w_{i+1}}$, and $S = N\cap S_{i}$ and
let \[{\mb Q} = {\mb P_w}|_N[-\dim_\C S_i]\cong 
\IC(X_w \cap N).\]  
The degree shift makes ${\mb Q}$ a perverse sheaf.

We give $N$ the stratification induced from $\{S_w\}$.  
Let $\gamma$ be a loop in $\wt\Lambda_S$ as in Theorem \ref{mainthm}.
We will show that the local system $\cM_S({\mb Q})$
has trivial monodromy around  
$\gamma$, or equivalently, that
$\cM_i$ has trivial monodromy around $\rho \circ \gamma$, 
where $\rho\colon \Lambda_N\to \Lambda$ is the inverse to the
homeomorphism given by restricting the natural projection
$T^*X|_N \to T^*N$ to $\Lambda\cap T^*X|_N$.
We can then apply Theorem \ref{mainthm} to show
that $\cM_{\{0\}}({\mb Q}) = \cM_{S_{i+1}}({\mb P}) \ne 0$.
 
There is an action of the torus $T=(\C^*)^{n+1}$ on $X$ which 
preserves the Schubert stratification, defined by acting on
$\C^{2n}$ via
\[(z_0, z_1, \dots, z_n)\cdot(x_1, \dots, x_{2n}) = 
(z_1^{-1}x_1, \dots, z_n^{-1}x_n, z_0z_nx_{n+1}, \dots, z_0z_1x_{2n}).\]
The $E_w$ are all fixed points of this action, the 
normal slice $N$ is preserved, and 
the induced action on $N$ is linear.  This action makes all Schubert variety
singularities conical, as the following lemma shows.

\begin{lemma}
For any word $w \in \{\alpha,\beta\}^n$, there is a 
homomorphism $\chi_w\colon \C^* \to T$ for which the induced 
action of $\C^*$ on $N_w$ is the conical action.
\end{lemma}
\begin{proof}
Let $\chi_w(z) = (z, z^{a_1}, \dots, z^{a_n})$, where $a_i = -1$ if 
$w(i) = \beta$, and $a_i = 0$ if $w(i) = \alpha$.  Checking that
this gives the required action is an exercise in local Grassmannian
coordinates, similar to arguments in \cite{BF} and \cite{B}.
\end{proof}

\begin{prop} The loop $\rho \circ \gamma$ 
is homotopic (in $\pi_1(\wt{\Lambda}_i)$)
to a loop generated by the action of a loop in $(\C^*)^{n+1}$.
\end{prop}
\begin{proof} 
Acting on $\gamma$ by the loop $\theta \mapsto \chi_{w_{i+1}}(e^{-i\theta})$,
we obtain the loop $\chi^{-1}_{w_{i+1}}\cdot \gamma(\theta) = (x, e^{i\theta}\xi)$; i.e.\ 
the point in $N$ stays fixed and the covector goes around by a conical action.  
We can then slide the loop $\rho \circ (\chi^{-1}_{w_i}\cdot \gamma)$ in $\wt\Lambda_i$
to get a loop of the form $(E_i, e^{i\theta}\xi')$ for some $\xi' \in T_{E_i}X$.
Acting now by $\chi_{w_i}(e^{i\theta})$ now produces a trivial loop.
\end{proof}

The fact that the monodromy of $\cM_i$ around $\rho \circ \gamma$
is trivial now follows from the following proposition.

\begin{prop} Let $X$ be a complex variety with an 
 action of a connected algebraic group $G$ and a $G$-invariant
stratification.  If ${\mb P} = \IC(\overline S)$ for some stratum $S$,
then the Morse local systems
$\cM_T({\mb P})$ are constant along loops generated by acting
by loops in $G$.
\end{prop}
\begin{proof} Since ${\mb P}$ is a $G$-equivariant perverse sheaf,
the local systems $\cM_T({\mb P})$ must also be $G$-equivariant.
\end{proof}

\end{document}